\documentclass[11pt]{article}

\usepackage{latexsym}
\usepackage{amsfonts}
\usepackage{amssymb}
\usepackage{amsmath}
\usepackage{mathrsfs}
\usepackage[all]{xy}
\usepackage{hyperref}

\newcommand{\bldd}{\mathbf{D}}
\newcommand{\bldr}{\mathbf{r}}
\newcommand{\blds}{\mathbf{s}}
\newcommand{\bldt}{\mathbf{t}}
\newcommand{\bldm}{\mathbf{m}}
\newcommand{\bldn}{\mathbf{n}}

\newcommand{\sgn}{\mbox{sign\,}}

\newcommand{\be}{\begin{equation}}
\newcommand{\ee}{\end{equation}}
\newcommand{\bea}{\begin{eqnarray}}
\newcommand{\eea}{\end{eqnarray}}
\newcommand{\bean}{\begin{eqnarray*}}
\newcommand{\eean}{\end{eqnarray*}}
\newcommand{\brray}{\begin{array}}
\newcommand{\erray}{\end{array}}

\newcommand{\newsection}[1]{\setcounter{equation}{0}
\setcounter{dfn}{0}
\section{#1}}

\newtheorem{dfn}{Definition}[section]
\newtheorem{thm}[dfn]{Theorem}
\newtheorem{lmma}[dfn]{Lemma}
\newtheorem{ppsn}[dfn]{Proposition}
\newtheorem{crlre}[dfn]{Corollary}
\newtheorem{xmpl}[dfn]{Example}
\newtheorem{rmrk}[dfn]{Remark}

\newcommand{\bdfn}{\begin{dfn}\rm}
\newcommand{\bthm}{\begin{thm}}
\newcommand{\blmma}{\begin{lmma}}
\newcommand{\bppsn}{\begin{ppsn}}
\newcommand{\bcrlre}{\begin{crlre}}
\newcommand{\bxmpl}{\begin{xmpl}}
\newcommand{\brmrk}{\begin{rmrk}\rm}

\newcommand{\edfn}{\end{dfn}}
\newcommand{\ethm}{\end{thm}}
\newcommand{\elmma}{\end{lmma}}
\newcommand{\eppsn}{\end{ppsn}}
\newcommand{\ecrlre}{\end{crlre}}
\newcommand{\exmpl}{\end{xmpl}}
\newcommand{\ermrk}{\end{rmrk}}

\newcommand{\bbc}{\mathbb{C}}
\newcommand{\bbz}{\mathbb{Z}}
\newcommand{\bbn}{\mathbb{N}}

\newcommand{\scrf}{\mathscr{F}}

\newcommand{\cla}{\mathcal{A}}

\newcommand{\clh}{\mathcal{H}}

\newcommand{\clk}{\mathcal{K}}
\newcommand{\cll}{\mathcal{L}}

\newcommand{\clg}{\mathcal{G}}

\newcommand{\prf}{\noindent{\it Proof\/}: }

\newcommand{\one}{{1\!\!1}}

\newcommand{\id}{\mbox{id}}

\def \qed { \mbox{}\hfill
$\Box$\vspace{1ex}}

\newcommand{\half}{\frac{1}{2}}

\begin{document}

%%%%%%%%%%%%%%%%%%%%%%%%%%%%%%%%%

%%%%%%%%%%%%%%%%%%%%%%%%%%%%%%%%%
\author{{\sc Partha Sarathi Chakraborty} and
{\sc Arupkumar Pal}}
\title{Characterization of $SU_q(\ell+1)$-equivariant spectral triples
for the odd dimensional
quantum spheres}
\maketitle

\begin{center}
\textit{Dedicated to Prof.\ K. R. Parthasarathy on his
seventieth birthday.}
\end{center}

%%%%%%%%%%%%%%%%%%%%%%%%%%%%%%%%%%
%%%%%  ABSTRACT
%%%%%%%%%%%%%%%%%%%%%%%%%%%%%%%%%%
 \begin{abstract}
The quantum group $SU_q(\ell+1)$ has a canonical action
on the odd dimensional sphere $S_q^{2\ell+1}$.
All odd spectral triples acting on the $L_2$ space of $S_q^{2\ell+1}$
and equivariant under this action have been characterized.
This characterization then leads to the construction of
an optimum family of equivariant spectral triples having
nontrivial $K$-homology class. These generalize the results
of Chakraborty \& Pal for $SU_q(2)$.
 \end{abstract}
{\bf AMS Subject Classification No.:} {\large 58}B{\large 34}, {\large
46}L{\large 87}, {\large
  19}K{\large 33}

%%%%%%%%%%%%%%%%%%%%%%%%%%%%%%%%%%%%%%%%%%%%%%%%%%%

% \tableofcontents

%%%%%%%%%%%%%%%%%%%%%%%%%%%%%%%%%%%%%%%%%%%%%%%%%%%
\newsection{Introduction}
Noncommutative differential geometry, which more commonly just
goes by the name noncommutative geometry, is an extension of
noncommutative topology and was initially developed in order
to handle certain spaces like the leaf space
of foliations or duals of groups whose topology or geometry
are difficult to study using machinery available in
classical geomtry or topology. As the subject developed,
more and more examples were found that are further away
from classical spaces but can be handled by noncommutative
geometric methods.

For quite sometime though, it was commonly believed that
quantum groups or their homogeneous spaces, which are
rather far removed from classical manifolds, are not
covered by the formalism of noncommutative geometry.
This notion changed with~(\cite{c-p1}), where the authors 
treated the case of the quantum $SU(2)$
group and found a   family of
spectral triples acting on its $L_2$-space
that are equivariant with respect to its natural (co)action.
This family is optimal, in the sense
that given any nontrivial equivariant Dirac operator $D$ acting
on the $L_2$ space, there exists a Dirac  operator $\widetilde{D}$
belonging to this family such that
$\mbox{sign\,}D$ is a compact perturbation of $\mbox{sign\,}\widetilde{D}$
and there exist reals $a$ and  $b$ such that
\[
|D| \leq a + b|\widetilde{D}|.
\]
Later Dabrowski et al~\cite{d-l-s-s-v} constructed another equivariant
spectral triple for $SU_q(2)$ on two copies of the $L_2$ space,
which was shown in \cite{c-p3} to be equivalent to a direct sum two 
spectral triples constructed in \cite{c-p1}.
Equivariant triples for the two dimensional Podles spheres
were constructed in \cite{d-s} and \cite{d-l-p-s}.
In a more recent paper (\cite{d-d-l}), D'Andrea et al gave a construction
of an equivariant spectral triple for the quantum four dimensional
spheres.

Our aim in the present paper is to look for higher
dimensional counterparts of the spectral triples
found in \cite{c-p1}.
We first formulate precisely what
one means by an equivariant spectral triple in a general
set up.
We then use a combinatorial method,
implicitly used in \cite{c-p1} and \cite{c-p2}, to 
characterize completely all odd spectral triples
acting on the $L_2$ space of the odd dimensional sphere $S_q^{2\ell+1}$
(see section~3 for the description)
and equivariant under the action of the $SU_q(\ell+1)$ group
\textbf{for all $\ell>1$}. This also leads to the construction
of an optimum family of equivariant nontrivial $(2\ell+1)$-summable
odd spectral triples
sharing all the properties of the triples for $SU_q(2)$ in \cite{c-p1}.

One should mention in this context that the 
construction by Hawkins \& Landi~(\cite{h-l}) does not deal 
with equivariance, and more importantly,
   they produce a (bounded) Fredholm
module, not a spectral triple, which is essential for
determining the smooth structure, giving a metric on the state space
and also help in computing the index map through
a local Chern character.

The paper is organised as follows.
In the next section, we will describe the
combinatorial method that was earlier used implicitly
in \cite{c-p1} and \cite{c-p2}.
We also formulate the notion of equivariance.
This has been done using the quantum group at
the function algebra level rather than passing on
to the quantum universal envelopping algebra level.
In section~3, we describe the $C^*$-algebra of continuous
functions on the odd dimensional quantum spheres
and state some of their relevant properties.
In section~4, we briefly recall the quantum group
$SU_q(\ell+1)$ and its representation theory.
In particular, we describe a nice basis for
the $L_2$ space and study the Clebsch-Gordon
coefficients. These are then used to describe the
action by left multiplication on the $L_2$ space explicitly.
In section~5, we give a description of the $L_2$ space of the sphere
and give a natural covriant representation on it.
In the last section, we give a precise characterization of
the singular values and of the sign, which helps us to
produce an optimal family of equivariant Dirac operators,
extending the results of \cite{c-p1} in the present case.

%%%%%%%%%%%%%%%%%%%%%%%%%%%%%%%%%%%%%%%%%%%%%%%%%%%
\newsection{Preliminaries}
%%%%%%%%%%%%%%%%%%%%%%%%%%%%%%%%%%%%%%%%%%%%%%%%%%%

%%%%%%%%%%%%%%%%%%%%%%%%%%%%%%%%%%%%%%%%%%%%%%%%%%%%%%%%%
\subsection{Equivariance}
%%%%%%%%%%%%%%%%%%%%%%%%%%%%%%%%%%%%%%%%%%%%%%%%%%%
Suppose $G$ is a compact group, quantum or classical,
and $\cla$ is a unital $C^*$-algebra. Assume that
$G$ has an action on $\cla$ given by
$\tau:\cla\rightarrow\cla\otimes C(G)$,
so that $(\id\otimes\Delta)\tau=(\tau\otimes\id)\tau$,
$\Delta$ being the coproduct.
In other words, we have a $C^*$-dynamical system $(\cla,G,\tau)$.

%%%%%%%%%%%%%%%%%%%%%%%%%%%%%%%%%%%%%%%%%%%%%
\bdfn
A covariant representation $(\pi,u)$
of $(\cla,G,\tau)$ consists of
a unital *-representation $\pi:\cla\rightarrow\cll(\clh)$,
a unitary representation $u$ of $G$ on $\clh$, i.e.\
a unitary element of the multiplier algebra $M(\clk(\clh)\otimes C(G))$
such that they obey the condition
$(\pi\otimes\id)\tau(a)=u(\pi(a)\otimes I)u^*$ for all $a\in\cla$.
\edfn
%%%%%%%%%%%%%%%%%%%%%%%%%%%%%%%%%%%%%%%%%%%%%

%%%%%%%%%%%%%%%%%%%%%%%%%%%%%%%%%%%%%%%%%%%%%
\bdfn
Suppose $(\cla, G,\tau)$ is a $C^*$-dynamical system.
An \textbf{odd $G$-equivariant spectral data 
for $(\cla, G,\tau)$} is
a quadruple $(\pi,u,\clh,D)$ where
\begin{enumerate}
\item
$(\pi,u)$ is a covariant representation of $(\cla, G,\tau)$ on 
 $\clh$,
\item
$\pi$  faithful,
\item
$u(D\otimes I)u^* =D\otimes I$,
\item
$(\pi,\clh,D)$ is an odd spectral triple.
\end{enumerate}
\edfn
%%%%%%%%%%%%%%%%%%%%%%%%%%%%%%%%%%%%%%%%%%%%%

%%%%%%%%%%%%%%%%%%%%%%%%%%%%%%%%%%%%%%%%%%%%%%%%%%%
\subsection{The general scheme}
%%%%%%%%%%%%%%%%%%%%%%%%%%%%%%%%%%%%%%%%%%%%%%%%%%%
Let $\clg$ be a graph and $(V_1,V_2)$ be a partition of the
vertex set. We say that $(V_1,V_2)$ admits an \textbf{infinite ladder}
if there exist infinite number of disjoint
paths each going from a point in $V_1$
to a point in $V_2$. Here two paths are disjoint means
that the set of vertices of one does not
intersect the set of vertices
of the other.

Suppose $\clh$ is a Hilbert space,
and $D$ is a self-adjoint operator on $\clh$ with compact resolvent.
Then $D$ admits a spectral resolution $\sum_{\gamma\in\Gamma} d_\gamma
P_\gamma$, where the $d_\gamma$'s are all distinct and each $P_\gamma$
is a finite dimensional projection.  Assume now onward that all the
$d_\gamma$'s are nonzero. Let $c$ be a positive real.  Let
us  define a graph $\clg_c$ as follows: take the vertex set $V$ to
be $\Gamma$.  Connect two vertices $\gamma$ and $\gamma'$ by an edge
if $|d_\gamma-d_{\gamma'}|<c$. Let $V^+=\{\gamma\in V: d_\gamma>0\}$ and
$V^-=\{\gamma\in V: d_\gamma<0\}$.  This will give us a partition of
$V$.
This partition has the important property that
$(V^+,V^-)$ does not admit an infinite ladder.
This is easy to see, because
if there is a path from $\gamma$ to $\delta$
and $d_\gamma>0$, $d_\delta<0$, then for some $\alpha$ on the path,
one must have $d_\alpha\in[-c,c]$.
Since the paths are disjoint, it
would contradict the compact resolvent condition.
We will call such a partition a \textbf{sign-determining} partition.

We will use this knowledge about the graph.
We start with a self-adjoint operator with
discrete spectrum. First choose a basis that
diagonalizes the operator $D$.
Next we use the action of the algebra elements on the basis
elements of $\clh$ and the boundedness of their
commutators with $D$.
This gives certain growth restrictions
on the $d_\gamma$'s. These will give us some information
about the edges in the graph. We exploit this knowledge
to characterize those partitions $(V_1,V_2)$ of the vertex set
that are sign-determining, i.\ e.\ do not admit any infinite ladder.
The sign of the operator $D$ must be of the form
$\sum_{\gamma\in V_1}P_\gamma-\sum_{\gamma\in V_2}P_\gamma$
where  $(V_1,V_2)$ is a sign-determining partition.
  Of course, for a given $c$, the graph
$\clg_c$ may have no edges, or too few edges (if the singular values
of $D$ happen to grow too fast), in which case, we will
be left with too many sign-determining
partitions.
Therefore, we would like to characterize
those partitions that are sign-determining
for all sufficiently large values of $c$.

In general the scheme outlined above will
be extremely difficult to carry out,
as the action of the algebra elements with
respect to the basis that diagonalizes $D$
may be quite complicated, and therefore
using boundedness of commutator conditions
will in general be very difficult. This is where equivariance
plays an extremely crucial role. It gives us a \textit{nice basis}
that diagonalizes $D$, so that the boundedness
of commutator conditions are simpler
and the subsequent steps become much more tractable.

%%%%%%%%%%%%%%%%%%%%%%%%%%%%%%%%%%%%%%%%%%%%
\section{The odd dimensional quantum spheres}
%%%%%%%%%%%%%%%%%%%%%%%%%%%%%%%%%%%%%%%%%%%%
%%%%%%%%%%%%%%%%%%%%%%%%%%%%%%%%%%%%%%%%%%%%
Let $q\in[0,1]$.
The $C^*$-algebra $A_\ell\equiv C(S_q^{2\ell+1})$ of the quantum
sphere $S_q^{2\ell+1}$
is the universal $C^*$-algebra generated by
elements
$z_1, z_2,\ldots, z_{\ell+1}$
satisfying the following relations (see~\cite{h-s}):
\bean
z_i z_j & =& qz_j z_i,\qquad 1\leq j<i\leq \ell+1,\\
z_i z_j^* & =& q z_j^* z_i ,\qquad 1\leq i\neq j\leq \ell+1,\\
z_i z_i^* - z_i^* z_i +
(1-q^{2})\sum_{k>i} z_k z_k^* &=& 0,\qquad \hspace{2em}1\leq i\leq \ell+1,\\
\sum_{i=1}^{\ell+1} z_i z_i^* &=& 1.
\eean
The $K$-theory groups for these algebras were
computed in  \cite{v-s} and \cite{h-s}.
\bppsn[\cite{v-s},\cite{h-s}]\label{ktheory}
$K_0(A_\ell)=K_1(A_\ell)=\bbz$.
\eppsn

The group $SU_q(\ell+1)$ has an action on $S_q^{2\ell+1}$.
Before we describe the action, let us recall the
definition of the quantum group $SU_q(\ell+1)$.
The $C^*$-algebra $C(SU_q(\ell+1))$ is the universal $C^*$-algebra
generated by  $\{ u_{ij} : i,j=1, \cdots ,\ell+1 \} $ 
obeying the relations:
\[ 
\sum_k u_{ki}^* u_{kj}= \delta_{ij} I,\quad
\sum_k u_{ik} u_{jk}^*= \delta_{ij} I
\] 
\[
\sum_{\mbox{\scriptsize{$k_i$'s distinct}}} 
\hspace{-1em}(-q)^{I(k_1,k_2, \cdots, k_{\ell+1})} 
   u_{j_1k_1}
\cdots u_{j_{\ell+1} k_{\ell+1}} = \begin{cases}
%          0 & \mbox{if $j_r=j_s$ for}\cr
% & \mbox{some $r \ne s$}\cr
      (-q)^{I( j_1,j_2, \cdots, j_{\ell+1})} & \mbox{$j_i$'s distinct}\cr
0 & \mbox{otherwise}
             \end{cases}
\]
where $I( k_1,k_2, \cdots, k_{\ell+1})$ is the number of inversions
in  $( k_1,k_2, \cdots, k_{\ell+1})$.
The group laws are given by the folowing maps:
\bean
\Delta(u_{ij})&=&\sum_k u_{ik}\otimes u_{kj}\qquad\mbox{ (Comultiplication)}\\
 S(u_{ij})&=&u_{ji}^*\qquad\mbox{ (Antipode)}\\
  \epsilon(u_{ij})&=&\delta_{ij}\qquad\mbox{ (Counit)}
\eean

The map
\[
\tau(z_i)=\sum_k z_k\otimes u_{ki}^*
\]
extends to a *-homomorphism $\tau$ from $A_\ell$
into $A_\ell\otimes C(SU_q(\ell+1))$ and
obeys $(\mbox{id}\otimes\Delta)\tau = (\tau\otimes\mbox{id})\tau$.
In other words this gives an action of $SU_q(\ell+1)$ on $A_\ell$.

%%%%%%%%%%%%%%%%%%%%%%%%%%%%%%%%%%%%%%%%%%%%%%%%%
\newsection{Preliminaries on $SU_q(\ell+1)$}
%%%%%%%%%%%%%%%%%%%%%%%%%%%%%%%%%%%%%%%%%%%%%%%%%
Our next job will be to get a description of the covariant
representation of the system $(A_\ell,SU_q(\ell+1),\tau)$
on $L_2(S_q^{2\ell+1})$. For this we need a few facts on
the representation theory of $SU_q(\ell+1)$. In the first
subsection we describe an important indexing of the basis
elements of the representation space of the irreducibles.
Then we describe the Clebsch-Gordon coefficients and compute
certain estimates. In the last subsection, we write down
explicitly the left multiplication operator on $L_2(SU_q(\ell+1))$.

%%%%%%%%%%%%%%%%%%%%%%%%%%%%%%%%%%%%%%%%%%%%%%%%%%%
\subsection{Gelfand-Tsetlin tableaux}
%%%%%%%%%%%%%%%%%%%%%%%%%%%%%%%%%%%%%%%%%%%%%%%%%%%
Irreducible unitary representations of the group
$SU_q(\ell+1)$ are indexed by
Young tableaux $\lambda=(\lambda_1,\ldots,\lambda_{\ell+1})$,
where $\lambda_i$'s are nonnegative integers,
$\lambda_1\geq \lambda_2\geq \ldots\geq \lambda_{\ell+1}$
(Theorem~1.5, \cite{w}).
Write $\clh_\lambda$ for the Hilbert space where
the irreducible $\lambda$ acts.
There are various ways of indexing the basis elements
of $\clh_\lambda$. The one we will use is due to Gelfand
and Tsetlin.
According to their prescription, basis elements for
$\clh_\lambda$ are parametrized by arrays of the form
\[
\bldr=\left(\begin{matrix}r_{11}&r_{12} &\cdots&r_{1,\ell}&r_{1,\ell+1}\cr
                      r_{21}&r_{22}&\cdots &r_{2,\ell}&\cr
                          &\cdots&&&\cr
                      r_{\ell,1}&r_{\ell,2}&&&\cr
                      r_{\ell+1,1}&&&&
                   \end{matrix}\right),
\]
where $r_{ij}$'s are integers satisfying
$r_{1j}=\lambda_j$ for $j=1,\ldots,\ell+1$,
$r_{ij}\geq r_{i+1,j}\geq r_{i,j+1}\geq 0$ for all $i$, $j$.
Such arrays are known as Gelfand-Tsetlin tableaux, to be abreviated
as GT tableaux for the rest of this section.
For a GT tableaux $\bldr$, the symbol $\bldr_{i\cdot}$ will denote its
$i$\raisebox{.4ex}{th} row.
It is well-known that two representations indexed respectively
by $\lambda$ and $\lambda'$ are equivalent if and only if
$\lambda_j-\lambda_j^\prime$ is independent of $j$ (\cite{w}).
Thus one gets an equivalence relation on the set of Young tableaux
$\{ \lambda=(\lambda_1,\ldots,\lambda_{\ell+1}):
\lambda_1\geq \lambda_2\geq \ldots\geq \lambda_{\ell+1}, \lambda_j\in\bbn\}$.
This, in turn, induces an equivalence relation on the set of
all GT tableaux $\Gamma=\{\bldr: r_{ij}\in\bbn,
  r_{ij}\geq r_{i+1,j}\geq r_{i,j+1}\}$: one says $\bldr$ and $\blds$
are equivalent if $r_{ij}-s_{ij}$ is independent of $i$ and $j$.
By $\Gamma$ we will mean the above set modulo this equivalence.

We will denote by $u^\lambda$ the irreducible unitary indexed by $\lambda$,
$\{e(\lambda,\bldr):\bldr_{1\cdot}=\lambda\}$ will denote an orthonormal basis
for $\clh_\lambda$ and $u^\lambda_{\bldr\blds}$ will stand for the matrix entries
of $u^\lambda$ in this basis. The symbol $\one$ will denote the Young tableaux
$(1,0,\ldots,0)$. We will often omit the symbol $\one$
and just write $u$ in order to denote $u^\one$.
Notice that any GT tableaux $\bldr$ with first row $\one$
must be, for some $i\in\{1,2,\ldots,\ell+1\}$, of the form $(r_{ab})$, where
\[
r_{ab}=\begin{cases}1 & \mbox{if $1\leq a\leq i$ and $b=1$},\cr
       0   &\mbox{otherwise.}\end{cases}
\]
Thus such a GT tableaux is uniquely determined by the integer $i$.
We will write just $i$ for this GT tableaux $\bldr$.
Thus for example, a typical matrix entry of $u^\one$ will be
written simply as $u_{ij}$.

Let $\bldr=(r_{ab})$ be a GT tableaux.
Let
$H_{ab}(\bldr):=r_{a+1,b}-r_{a,b+1}$ and
$V_{ab}(\bldr):=r_{ab}-r_{a+1,b}$.
An element $\bldr$ of $\Gamma$ is completely
specified by the following differences
\[
\bldd(\bldr)=\left(\begin{matrix}
                    V_{11}(\bldr)&H_{11}(\bldr)
         &H_{12}(\bldr)&\cdots&H_{1,\ell-1}(\bldr)&H_{1,\ell}(\bldr)\cr
  V_{21}(\bldr)&H_{21}(\bldr)&H_{22}(\bldr)&\cdots&H_{2,\ell-1}(\bldr)&\cr
         &\cdots&&&&\cr
       V_{\ell,1}(\bldr)&H_{\ell,1}(\bldr)&&&&
    \end{matrix}\right).
\]
The differences satisfy the following inequalities
\be\label{ineq}
\sum_{k=0}^b H_{a-k,k+1}(\bldr)\leq V_{a+1,1}(\bldr)
       +\sum_{k=0}^b H_{a-k+1,k+1}(\bldr),\quad
1\leq a\leq \ell,\;\;0\leq b\leq a-1.
\ee
Conversely, if one has an array of the form
\[
\left(\begin{matrix}V_{11}&H_{11}&H_{12}&\cdots&H_{1,\ell-1}&H_{1,\ell}\cr
       V_{21}&H_{21}&H_{22}&\cdots&H_{2,\ell-1}&\cr
         &\cdots&&&&\cr
       V_{\ell,1}&H_{\ell,1}&&&&
    \end{matrix}\right),
\]
where $V_{ij}$'s and $H_{ij}$'s are in $\bbn$ and obey
the inequalities~(\ref{ineq}), then the above array is of the form
$\bldd(\bldr)$ for some GT tableaux $\bldr$. Thus the quantities
$V_{a1}$ and $H_{ab}$ give a coordinate system for elements in $\Gamma$.
The following diagram explains this new coordinate system.
The hollow circles stand for the $r_{ij}$'s.
The entries are decreasing along the direction of the arrows,
and the $V_{ij}$'s and the $H_{ij}$'s are the difference
between the two endpoints of the corresponding arrows.

%%%%%%%%%%%%%%%%%%%%%%%%%%%%%%%%%%%%%%%%%%%%%%%%
\hspace*{100pt}
  \def\labelstyle{\scriptstyle}
 
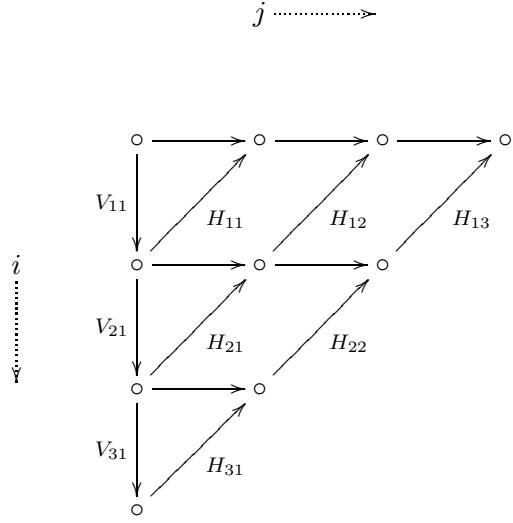
\xymatrix@C=35pt@R=35pt{
   &  &  j\ar@{.>}[r] &&\\
   & \circ\ar@{->}[r]\ar@{->}[d]_{V_{11}} &  \circ\ar@{->}[r]
                       & \circ\ar@{->}[r] &\circ\\
i\ar@{.>}[d] &  \circ\ar@{->}[r]\ar@{->}[d]_{V_{21}}\ar@{->}[ur]_{H_{11}} &
    \circ\ar@{->}[r]\ar@{->}[ur]_{H_{12}} & \circ\ar@{->}[ur]_{H_{13}} & \\
   & \circ\ar@{->}[r]\ar@{->}[d]_{V_{31}}\ar@{->}[ur]_{H_{21}} &
               \circ\ar@{->}[ur]_{H_{22}} &\\
   &  \circ\ar@{->}[ur]_{H_{31}} & }\\
%%%%%%%%%%%%%%%%%%%%%%%%%%%%%%%%%%%%%%%%%%%%%%%%

%%%%%%%%%%%%%%%%%%%%%%%%%%%%%%%%%%%%%%%%%%%%%%%%%%%%%%%%%%%
%% CG coefficients
%%%%%%%%%%%%%%%%%%%%%%%%%%%%%%%%%%%%%%%%%%%%%%%%%%%%%%%%%%%
\subsection{Clebsch-Gordon coefficients}
%%%%%%%%%%%%%%%%%%%%%%%%%%%%%%%%%%%%%%%%%%%%%%%%%%%%%%%%%%%
In this subsection, we recall the Clebsch-Gordon coefficients
for the group $SU_q(\ell+1)$. This will be important in writing
down the natural representation of $C(S_q^{2\ell+1})$ on $L_2(S_q^{2\ell+1})$ explicitly.

Look at the representation $u^\one\otimes u^\lambda$
acting on $\clh_\one\otimes\clh_\lambda$.
The representation decomposes as a direct sum
$\oplus_\mu u^\mu$, i.e.\ one has a corresponding
decomposition $\oplus_\mu\clh_\mu$ of $\clh_\one\otimes\clh_\lambda$.
Thus one has two orthonormal bases
$\{e^\mu_\blds\}$ and $\{e^\one_i\otimes e^\lambda_\bldr\}$.
The Clebsch-Gordon coefficient $C_q(\one,\lambda,\mu;i,\bldr,\blds)$
is defined to be the inner product
$\langle e^\mu_\blds, e^\one_i\otimes e^\lambda_\bldr\rangle$.
Since $\one$, $\lambda$ and $\mu$ are just the first rows of
$i$, $\bldr$ and $\blds$ respectively, we will often denote
the above quantity just by $C_q(i,\bldr,\blds)$.

Next, we will compute the quantities $C_q(i,\bldr,\blds)$.  We will
use the calculations given in (\cite{k-s}, pp.\ 220), keeping in mind
that for our case (i.e.\ for $SU_q(\ell+1)$), the top right entry of
the GT tableaux is zero.

Let $M=(m_1,m_2,\ldots,m_i)\in\bbn^i$ be such that $1\leq m_j\leq \ell+2-j$.
Denote by $M(\bldr)$ the tableaux $\blds$ defined by
\be\label{movenotation}
s_{jk}=\begin{cases}r_{jk}+1 & \mbox{if $k=m_j$, $1\leq j\leq i$},\cr
         r_{jk} & \mbox{otherwise.}
 \end{cases}
\ee
With this notation, observe now that
$C_q(i,\bldr,\blds)$ will be zero unless $\blds$ is
$M(\bldr)$ for some $M\in\bbn^i$.
(One has to keep in mind though that not all tableaux of the form $M(\bldr)$
is a valid GT tableaux)

From (\cite{k-s}, pp.\ 220), we have
\be\label{cgc1}
C_q(i,\bldr,M(\bldr))=\prod_{a=1}^{i-1}
\left\langle \brray{ll}
                (1,\mathbf{0}) &\bldr_{a\cdot} \cr
                (1,\mathbf{0}) &\bldr_{a+1\cdot}
             \erray\left| \brray{l}
                   \bldr_{a\cdot}+e_{m_a}\cr
                   \bldr_{a+1\cdot}+e_{m_{a+1}}
              \erray\right.\right\rangle
\times
\left\langle \brray{ll}
                (1,\mathbf{0}) &\bldr_{i\cdot} \cr
                (0,\mathbf{0}) &\bldr_{i+1\cdot}
             \erray\left| \brray{l}
                   \bldr_{i\cdot}+e_{m_i}\cr
                   \bldr_{i+1\cdot}
              \erray\right.\right\rangle,
\ee
where $e_k$ stands for a vector (in the appropriate space) whose
$k$\raisebox{.4ex}{th} coordinate is 1 and the rest are all zero, 
$\bldr_{j\cdot}$ stands for the j\raisebox{.4ex}{th} row of the
tableaux $\bldr$,
and
\bea
\left\langle \brray{ll}
                (1,\mathbf{0}) &\bldr_{a\cdot} \cr
                (1,\mathbf{0}) &\bldr_{a+1\cdot}
             \erray\left| \brray{l}
                   \bldr_{a\cdot}+e_j\cr
                   \bldr_{a+1\cdot}+e_k
              \erray\right.\right\rangle^2
&=&
q^{-r_{aj}+r_{a+1,k} - k+j}
\times
\prod_{{i=1}\atop{i\neq j}}^{\ell+2-a}
   \frac{[r_{a,i}-r_{a+1,k}-i+k]_q  }{[r_{a,i}-r_{a,j}-i+j]_q} \nonumber \\
&& \times
\prod_{{i=1}\atop{i\neq k}}^{\ell+1-a}
   \frac{[r_{a+1,i}-r_{a,j}-i+j-1]_q  }{[r_{a+1,i}-r_{a+1,k}-i+k-1]_q},\label{corrected_1}\\
\left\langle \brray{ll}
                (1,\mathbf{0}) &\bldr_{a\cdot} \cr
                (0,\mathbf{0}) &\bldr_{a+1\cdot}
             \erray\left| \brray{l}
                   \bldr_{a\cdot}+e_j\cr
                   \bldr_{a+1\cdot}
              \erray\right.\right\rangle^2
&=& q^{\left(1-j+\sum_{i=1}^{\ell+1-a}r_{a+1,i} -
           \sum_{{i=1}\atop{i\neq j}}^{\ell+2-a}r_{a,i}\right)} \nonumber \\
&&
\times \left(
\frac{\prod_{i=1}^{\ell+1-a}[r_{a+1,i}-r_{aj}-i+j-1]_q  }
{\prod_{{i=1}\atop{i\neq j}}^{\ell+2-a}[r_{a,i}-r_{aj}-i+j]_q  }\right),
\label{corrected_2}
\eea
where for an integer $n$, $[n]_q$ denotes the $q$-number $(q^n-q^{-n})/(q-q^{-1})$.
After some lengthy but straightforward computations,
we get the following two relations:
\be
\left|
\left\langle \brray{ll}
                (1,\mathbf{0}) &\bldr_{a\cdot} \cr
                (1,\mathbf{0}) &\bldr_{a+1\cdot}
             \erray\left| \brray{l}
                   \bldr_{a\cdot}+e_j\cr
                   \bldr_{a+1\cdot}+e_k
              \erray\right.\right\rangle
\right| = A'q^A,
\ee
\be
\left|
\left\langle \brray{ll}
                (1,\mathbf{0}) &\bldr_{a\cdot} \cr
                (0,\mathbf{0}) &\bldr_{a+1\cdot}
             \erray\left| \brray{l}
                   \bldr_{a\cdot}+e_j\cr
                   \bldr_{a+1\cdot}
              \erray\right.\right\rangle
\right| = B'q^B,
\ee
where
\bea
A&=&\begin{cases}\displaystyle{\sum_{j\wedge k < b < j\vee k}(r_{a+1,b}-r_{a,b})}
              +(r_{a+1,j\wedge k}-r_{a,j\vee k})  & \mbox{if $j\neq k$},\cr
   0 & \mbox{if $j=k$}.\end{cases} \cr
&=& \sum_{j\wedge k \leq b < j\vee k}(r_{a+1,b}-r_{a,b+1})
    +2 \sum_{k < b < j}(r_{a,b}-r_{a+1,b}) \cr
&=&  \sum_{j\wedge k \leq b < j\vee k}H_{ab}(\bldr)
             + 2 \sum_{k < b < j}V_{ab}(\bldr).\label{cgc2}\\
B &=&  \sum_{j \leq b < \ell+2-a}H_{ab}(\bldr),\label{cgc3}
\eea
and $A'$ and $B'$ both lie between two positive constants
independent of $\bldr$, $a$, $j$ and $k$
(Here and elsewhere in this paper, an empty summation
would always mean zero).

Combining these, one gets
\be \label{cgc4}
C_q(i,\bldr, M(\bldr))=P\cdot q^{C(i,\bldr,M)},
\ee
where
\be \label{cgc5}
C(i,\bldr,M)=\sum_{a=1}^{i-1}\left(
   \sum_{m_a\wedge m_{a+1} \leq b < m_a\vee m_{a+1}}H_{ab}(\bldr)
   +2 \sum_{m_{a+1} < b < m_a}V_{ab}(\bldr)\right)
+\sum_{m_i \leq b < \ell+2-i}H_{ib}(\bldr),
\ee
and $P$ lies between two positive constants
that are independent of $i$, $\bldr$ and $M$.

\brmrk
The formulae (\ref{corrected_1}) and (\ref{corrected_2})
are obtained from equations~(45) and (46), page 220, \cite{k-s}
by replacing $q$ with $q^{-1}$. Equation~(45) is a special
case of the more general formula (48), page 221, \cite{k-s}.
However, there is a small error in equation~(48) there.
The correct form can be found in equations~(3.1, 3.2a, 3.2b)
in \cite{a-s}. That correction has been incorporated in
equations~(\ref{corrected_1}) and (\ref{corrected_2}) here.
\ermrk

%%%%%%%%%%%%%%%%%%%%%%%%%%%%%%%%%%%%%%%%%%%%%%%%%%%%%%%%%%%%%%%%%%%
\subsection{Left multiplication operators}
%%%%%%%%%%%%%%%%%%%%%%%%%%%%%%%%%%%%%%%%%%%%%%%%%%%%%%%%%%%%%%%%%%%
We next write down the representation of $C(SU_q(\ell+1))$ on $L_2(SU_q(\ell+1))$ by 
left multiplication. Later we will work with a certain
restriction of this representation.

The matrix entries $u^\lambda_{\bldr\blds}$ form a complete orthogonal set
of vectors in $L_2(SU_q(\ell+1))$. Write $e^\lambda_{\bldr\blds}$ for
$\|u^\lambda_{\bldr\blds}\|^{-1}u^\lambda_{\bldr\blds}$.
Then the $e^\lambda_{\bldr\blds}$'s form a complete orthonormal basis
for $L_2(SU_q(\ell+1))$. Let $\pi$ denote the representation of $\cla$ on
$L_2(SU_q(\ell+1))$ by left multiplications. We will now derive an expression for
$\pi(u_{ij})e^\lambda_{\bldr\blds}$.

From the definition of matrix entries and that of the CG coefficients,
one gets
\be \label{cb1}
u^\rho e(\rho,\bldt)=\sum_\blds u^\rho_{\blds\bldt}e(\rho,\blds),
\ee
\be \label{cb2}
e(\mu,\bldn)=\sum_{j,\blds}C_q(j,\blds,\bldn)e(\one,j)\otimes e(\lambda,\blds).
\ee
Apply $u\otimes u^\lambda$ on both sides and note that
$u\otimes u^\lambda$ acts on $e(\mu,\bldn)$ as $u^\mu$:
\be \label{cb3}
\sum_\bldm u^\mu_{\bldm\bldn}e(\mu,\bldm)=
\sum_{j,\blds}\sum_{i,\bldr}C_q(j,\blds,\bldn)
 u_{ij}u^\lambda_{\bldr\blds}e(\one,i)\otimes e(\lambda,\bldr).
\ee
Next, use (\ref{cb2}) to expand $e(\mu,\bldm)$ on the left hand side to get
\be
\sum_{i,\bldr,\bldm} u^\mu_{\bldm\bldn}
   C_q(i,\bldr,\bldm)e(\one,i)\otimes e(\lambda,\bldr)
=
\sum_{j,\blds}\sum_{i,\bldr}C_q(j,\blds,\bldn)
        u_{ij}u^\lambda_{\bldr\blds}e(\one,i)\otimes
         e(\lambda,\bldr).
\ee
Equating coefficients, one gets
\be
\sum_{\bldm} C_q(i,\bldr,\bldm)u^\mu_{\bldm\bldn}
=
\sum_{j,\blds}C_q(j,\blds,\bldn)
        u_{ij}u^\lambda_{\bldr\blds}.
\ee
Now using orthogonality of the matrix
$(\!(C_q(\one,\lambda,\mu;j,\blds,\bldn))\!)_{(\mu,\bldn),(j,\blds)}$,
we obtain
\be\label{alg_left_mult}
u_{ij}u^\lambda_{\bldr\blds}
= \sum_{\mu,\bldm,\bldn}
  C_q(i,\bldr,\bldm)C_q(j,\blds,\bldn)u^\mu_{\bldm\bldn}.
\ee
From (\cite{k-s}, pp.\ 441), one has
$\|u^\lambda_{\bldr\blds}\|=d_\lambda^{-\half}q^{-\psi(\bldr)}$,
where
\[
\psi(\bldr)=-\frac{\ell}{2}\sum_{j=1}^{\ell+1}r_{1j}
   + \sum_{i=2}^{\ell+1}\sum_{j=1}^{\ell+2-i}r_{ij},
   \qquad
d_\lambda=\sum_{\bldr:\bldr_1=\lambda} q^{2\psi(\bldr)}
\]

Therefore
\be\label{left_mult}
\pi(u_{ij})e^\lambda_{\bldr\blds}
= \sum_{\mu,\bldm,\bldn}
  C_q(\one,\lambda,\mu;i,\bldr,\bldm)C_q(\one,\lambda,\mu;j,\blds,\bldn)
   d_\lambda^\half d_\mu^{-\half}q^{\psi(\bldr)-\psi(\bldm)}
   e^\mu_{\bldm\bldn}.
\ee

Write
\be
\kappa(\bldr,\bldm)=
 d_\lambda^\half d_\mu^{-\half}q^{\psi(\bldr)-\psi(\bldm)}.
\ee
%%%%%%%%%%%%%%%%%%%%%%%%%%%%%%%%%%%%%%%%%%%%%%%%%%%%%%%%%%%%%%%%%%%
\blmma\label{krmbound}
There exist constants $K_2>K_1>0$ such that
$K_1< \kappa(\bldr, M(\bldr))<K_2$ for all $\bldr$.
\elmma
%%%%%%%%%%%%%%%%%%%%%%%%%%%%%%%%%%%%%%%%%%%%%%%%%%%%%%%%%%%%%%%%%%%
\prf
Observe that
\[
\psi(\bldr)=(\rho,\lambda(\bldr))
=-\frac{\ell}{2}\sum_{j=1}^{\ell+1}r_{1j}+
   \sum_{i=2}^{\ell+1}\sum_{j=1}^{\ell+2-i}r_{ij}.
\]
Therefore
\[
\min\{\psi(\bldr):\bldr_1=\lambda\}
  = -\frac{\ell}{2}\sum_1^\ell\lambda_i + \sum_{k=2}^\ell(k-1)\lambda_k.
\]
This implies that
\[
d_\lambda^\half=q^{-\frac{\ell}{2}\sum_1^\ell\lambda_i + \sum_{k=2}^\ell(k-1)\lambda_k}(1+o(q)),
\]
which gives us
\[
\left(\frac{d_\lambda}{d_{\lambda+e_k}}\right)^\half=
q^{\frac{\ell}{2}-M_1+1}(1+o(q)).
\]
Next,
\[
q^{\psi(\bldr)-\psi(M(\bldr))}=
 q^{-\frac{\ell}{2}\sum_{j=1}^{\ell+1}r_{1j}+
   \sum_{i=2}^{\ell+1}\sum_{j=1}^{\ell+2-i}r_{ij}
    +\frac{\ell}{2}(\sum_{j=1}^{\ell+1}r_{1j}+1)
   - (\sum_{i=2}^{\ell+1}\sum_{j=1}^{\ell+2-i}r_{ij}+i-1) }
   =q^{\frac{\ell}{2}-i+1}.
\]
Thus
\[
\kappa(\bldr,M(\bldr))=q^{\ell-i-M_1+2}(1+o(q)).
\]
Hence the conclusion follows.\qed

%%%%%%%%%%%%%%%%%%%%%%%%%%%%%%%%%%%%%%%%%%%%%%%%%%%%%%%%%%%%%%
\newsection{Covariant representation}
%%%%%%%%%%%%%%%%%%%%%%%%%%%%%%%%%%%%%%%%%%%%%%%%%%%%%%%%%%%%%%
%%%%%%%%%%%%%%%%%%%%%%%%%%%%%%%%%%%%%%%%%%%%
%%%%%%%%%%%%%%%%%%%%%%%%%%%%%%%%%%%%%%%%%%%%
Let us write $G$ for $SU_q(\ell+1)$ and $H$ for $SU_q(\ell)$.
$H$ is a subgroup of $G$. This means 
that there is a $C^*$-epimorphism $\phi:C(G)\rightarrow C(H)$
  obeying $\Delta_H\phi=(\phi\otimes\phi)\Delta_G$.
In such a case, one defines the quotient space $G\backslash H$ by
\[
C(G\backslash H):=\{a\in C(G): (\phi\otimes id)\Delta (a)=I\otimes a\}.
\]
The group $G$ has a canonical right action
$C(G\verb1\1H)\rightarrow C(G\verb1\1H)\otimes C(G)$
coming from the restriction of the comultiplication $\Delta$
to $C(G\verb1\1H)$.
Let $\rho$ denote the restriction of the Haar state on $C(G)$
to $C(G\verb1\1H)$.
Then clearly one has
$(\rho\otimes id)\Delta (a) = \rho(a)I$,
which means $\rho$ is the invariant state for $C(G\verb1\1H)$.
This also means that $L_2(G\verb1\1H)=L_2(\rho)$ is just the
closure of $C(G\verb1\1H)$ in $L_2(G)$.
(For a formulation of quotient spaces etc.\ in the context
of compact quantum groups, see~\cite{po})

Now suppose we make the following explicit choice of $\phi$.
Let $u^\one$ denote the fundamental unitary for $G$,
i.\ e.\ the irreducible unitary representation corresponding to the
Young tableaux $\one=(1,0,\ldots,0)$.
Similarly write $v^\one$ for the fundamental unitary for $H$.
Fix some bases for the corresponding representation spaces.
Then $C(G)$ is the $C^*$-algebra generated by the matrix
entries $\{u^\one_{ij}\}$ and $C(H)$ is the
$C^*$-algebra generated by the matrix
entries $\{v^\one_{ij}\}$.
Now define $\phi$ by
\be
\phi(u^\one_{ij})=\begin{cases} I & \mbox{if $i=j=1$},\cr
        v^\one_{i-1,j-1} & \mbox{if $2\leq i,j\leq \ell+1$},\cr
        0 & \mbox{otherwise.}
       \end{cases}
\ee
Then $C(G\verb1\1H)$ is the $C^*$-subalgebra of
$C(G)$ generated by the entries $u_{1,j}$ for $1\leq j\leq \ell+1$.
Define $\theta:A_\ell\rightarrow C(G\verb1\1H)$
by 
\[
 \theta(z_i)=q^{-i+1}u^*_{1,i}.
\]
This gives an isomorphism between $C(G\verb1\1H)$ and $A_\ell$
and the following diagram commutes:\\[1ex]
\[
\def\labelstyle{\scriptstyle}
  \xymatrix@C=23pt@R=50pt{
 A_\ell\ar[d]_{\theta}\ar[r]_-{\tau} & A_\ell\otimes C(G)\ar[d]_{\theta\otimes\id}&\\
C(G\backslash H)\ar[r]_-{\Delta}  & C(G\backslash H)\otimes C(G)\\
}
\]
% (one recovers the relations for the generators of $C(S_q^{2\ell+1})$
% if one sets $z_i=q^{-i+1}u^*_{1,i})$.
In other words,
$(A_\ell, G,\tau)$ is the quotient space $G\backslash H$.
As we shall see shortly, this choice of $\phi$
will make $L_2(G\verb1\1H)$ 
a span of certain rows of the $e_{\bldr,\blds}$'s
and this in turn will help us make use of the calculations
already done in the initial sections.

\bppsn
Assume $\ell>1$.
The right regular representation $u$ of $G$ keeps
$L_2(G\verb1\1H)$ invariant, and the restriction of $u$ to
$L_2(G\verb1\1H)$ decomposes as a direct sum of exactly one copy
of each of the irreducibles given by the young tableaux
$\lambda_{n,k}:=(n+k, k,k,\ldots, k,0)$, with $n,k\in\bbn$.
\eppsn
\prf
Write $\sigma$ for the composition $h_H\circ\phi$
where $h_H$ is the Haar state for $H$.
From the description of $C(G\verb1\1H)$, it follows
that
\bean
C(G\verb1\1H) &=& \{a\in C(G): (\sigma \otimes id)\Delta (a)=a\}\\
& =&\{(\sigma \otimes id)\Delta (a): a\in C(G)\}.
\eean
Now
the map $a\mapsto \sigma\ast a:=(\sigma\otimes id)\Delta(a)$
on $C(G)$ extends to a bounded linear operator $L_\sigma$ on $L_2(G)$
(lemma~3.1, \cite{pa}), and it is easy to see that
$L_\sigma^2=L_\sigma$. It follows then that
$L_2(G\verb1\1H)=\ker(L_\sigma -I)=\mbox{ran}\,L_\sigma$.
From the discussion preceeding theorem~3.3, \cite{pa},
it now follows that $u$ keeps $L_2(G\verb1\1H)$ invariant and in fact
the restriction of $u$ to $L_2(G\verb1\1H)$ is the representation
induced by the trivial repersentation of $H$.
From the analogue of Frobenius reciprocity theorem for
compact quantum groups (theorem~3.3, \cite{pa}) it now
follows that the multiplicity of any irreducible $u^\lambda$
in it would be same as the multiplicity of the trivial
representation of $H$ in the restriction of $u^\lambda$ to $H$.
But from the representation theory of $SU_q(\ell+1)$,
we know that  the restriction of $u^\lambda$ to $SU_q(\ell)$
decomposes into a direct sum of one copy of
each irreducible $\mu:(\mu_1\geq \mu_2\geq \ldots \geq\mu_\ell)$
of $SU_q(\ell)$ for which
\be\label{induced}
\lambda_1\geq \mu_1 \geq \lambda_2\geq \mu_2\geq \ldots
         \geq\lambda_\ell \geq \mu_\ell \geq 0.
\ee
Now the trivial representation of $SU_q(\ell)$ is indexed
by Young tableaux of the form
$\mu:(k,k,\ldots,k)$ where $k\in\bbn$.
But such a $\mu$ will obey the restriction~\ref{induced} above
if and only if $\lambda$ is of the form
$(n+k,k,k,\ldots,k,0)$.
\qed
\brmrk
For the case $\ell=1$, the restriction of the irreducible
$(n,0)$ to the trivial subgroup decomposes into $n+1$ copies
of the trivial  representation. Therefore, in this case,
$L_2(S_q^3)$ decomposes into a direct sum of $n+1$ copies of
each representation $(n,0)$.
\ermrk

\bppsn
Let $\Gamma_0$ be the set of all GT tableaux $\bldr^{nk}$
given by
\[
r^{nk}_{ij}=\begin{cases} n+k & \mbox{if $i=j=1$},\cr
                0  & \mbox{if $i=1$, $j=\ell+1$},\cr
                k  & \mbox{otherwise},\end{cases}
\]
for some $n,k \in \bbn$.
Let $\Gamma_0^{nk}$ be the set of all GT tableaux with
top row $(n+k,k,\ldots,k,0)$.
Then the family of vectors
\[
\{e_{\bldr^{nk},\blds}: n,k\in\bbn,\, \blds\in\Gamma_0^{nk}\}
\]
form a complete
orthonormal basis for $L_2(G\verb1\1H)$.
\eppsn
\prf
Let $A$ be the linear span of the elements
$\{u_{\bldr^{n,k},\blds}: n,k\in\bbn, \blds\in\Gamma_0^{n,k}\}$.
Clearly  the closure of $A$ in $L_2(G)$ is the closed
linear span of $\{e_{\bldr^{nk},\blds}: n,k\in\bbn,\, \blds\in\Gamma_0^{nk}\}$.
It is also immdiate that the restriction of the
right regular representation to the above subspace
is a direct sum of one copy of each of the irreducibles
$(n+k,k,k,\ldots,k,0)$.

We will next show that for any $a\in A$, $u_{1j}a$ and $u_{1j}^*$ a are also
in $A$.
Take $a=u_{\bldr^{n,k},\blds}$. Use equation~(\ref{alg_left_mult}) to get
\bea
u_{1,j}u_{\bldr^{n,k},\blds} &=&
  \sum_{M, M'}
  C_q(1,\bldr^{n,k},M(\bldr^{n,k}))C_q(j,\blds,M'(\blds))
      u_{M(\bldr^{n,k}),M'(\blds)}\cr
 &=& \sum_{M'}
  C_q(1,\bldr^{n,k},M_{11}(\bldr^{n,k}))C_q(j,\blds,M'(\blds))
      u_{M_{11}(\bldr^{n,k}),M'(\blds)} \cr
&&   +
   \sum_{M''}
  C_q(1,\bldr^{n,k},M_{\ell+1,1}(\bldr^{n,k}))C_q(j,\blds,M''(\blds))
      u_{M_{\ell+1,1}(\bldr^{n,k}),M''(\blds)} \cr
&=&\sum_{M'}
  C_q(1,\bldr^{n,k},\bldr^{n+1,k})C_q(j,\blds,M'(\blds))
      u_{\bldr^{n+1,k},M'(\blds)}\cr
&&   +
   \sum_{M''}
  C_q(1,\bldr^{n,k},\bldr^{n,k-1}))C_q(j,\blds,M''(\blds))
      u_{\bldr^{n,k-1},M''(\blds)},
\eea
where the first sum is over all moves $M'\in\bbn^{j}$
whose first coordinate is 1 and the second
sum is over all moves $M''\in\bbn^{j}$
whose first coordinate is $\ell+1$.
Thus $u_{1j}a\in A$.

Next, note that if
$\langle u_{1j}^* e_{\bldr^{n,k},\blds}, e_{\bldr',\blds'}\rangle\neq 0$,
then one must have
$\bldr'=\bldr^{n-1,k}$ or $\bldr'=\bldr^{n,k+1}$.
Therefore it follows that $u_{1j}^* u_{\bldr^{n,k},\blds}$
is a linear combination of the $u_{\bldr^{n-1,k},\blds}$
$u_{\bldr^{n,k+1},\blds}$'s, and hence belongs to $A$.
Since $A$ contains the element $u_{\mathbf{0},\mathbf{0}}=1$,
it contains $u_{1j}$ and $u_{ij}^*$. Thus $A$ contains the $*$-algebra
$B$ generated by the $u_{1j}$'s.
But by the previous theorem, restriction of the right regular representation
to the $L_2$ closure $L_2(G\verb1\1H)$ of $B$ also decomposes
as a direct sum of one copy
of each of the irreducibles $(n+k,k,\ldots,k,0)$.
So it follows that $L_2(G\verb1\1H)$ is equal to
the subspace stated in the theorem.
\qed

Thus the right regular representation u restricts to the subspace
$L_2(G\backslash H)$ and it also follows from the above discussion
that the restriction of the left multiplication to $C(G\backslash H)$
keeps $L_2(G\backslash H)$ invariant.
Let us denote the restriction of $u$ to $L_2(G\backslash H)$
by $\hat{u}$ and the restriction of $\pi$ to $C(G\backslash H)$
viewed as a map on $L_2(G\backslash H)$ by $\hat{\pi}$.
It is easy to check that $(\hat{\pi},\hat{u})$ is a covariant representation
for the system $(A_\ell,G,\tau)$.

%%%%%%%%%%%%%%%%%%%%%%%%%%%%%%%%%%%%%%%%%%%%%%%%%%%%%%%%%%%%%%
\newsection{Equivariant spectral triples}
%%%%%%%%%%%%%%%%%%%%%%%%%%%%%%%%%%%%%%%%%%%%%%%%%%%%%%%%%%%%%%
The following lemma is straightforward.
\blmma
Let $D$ be a self-adjoint operator with compact resolvent on $L_2(G\verb1\1H)$
that is equivariant with respect to the covariant representation $(\hat{\pi},\hat{u})$
then it is of the form
\[
e_{\bldr,\blds}\mapsto d(\bldr)e_{\bldr,\blds},\quad \bldr\in\Gamma_0.
\]
\elmma

For such a $D$, one can then write down the commutatots
with algebra elements:
\be\label{bdd_comm}
[D,\pi(u_{ij})]e^\lambda_{\bldr\blds}=
\sum (d(\bldm)-d(\bldr))C_q(\one,\lambda,\mu;i,\bldr,\bldm)
    C_q(\one,\lambda,\mu;j,\blds,\bldn)
\kappa(\bldr,\bldm)e^\mu_{\bldm\bldn}.
\ee
Therefore the condition for boundedness of commutators reads
as follows:
\be \label{eqbdd1}
|(d(\bldm)-d(\bldr))C_q(\one,\lambda,\mu;i,\bldr,\bldm)
   C_q(\one,\lambda,\mu;j,\blds,\bldn)
\kappa(\bldr,\bldm)|<c,
\ee
where $c$ is independent of $i$, $j$, $\lambda$, $\mu$, $\bldr$, $\blds$, $\bldm$ and $\bldn$.

Using lemma~\ref{krmbound}, we get
\be\label{eqbdd2}
|(d(\bldm)-d(\bldr))C_q(\one,\lambda,\mu;i,\bldr,\bldm)
  C_q(\one,\lambda,\mu;j,\blds,\bldn)|<c.
\ee
Choosing $j$, $\blds$ and $\bldn$ suitably, one can ensure that
(\ref{eqbdd2}) implies the following:
\be\label{eqbdd3}
|(d(\bldm)-d(\bldr))C_q(\one,\lambda,\mu;i,\bldr,\bldm)|<c.
\ee
It follows from~(\ref{bdd_comm}) that this condition is also sufficient for
the boundedness of the commutators $[D, u_{ij}]$.

From (\ref{cgc4}), one gets
%%%%%%%%%%%%%%%%%%%%%%%%%%%%%%%%%%%%%%%%%%%%%%%%%%%%%
\be \label{eqbdd4}
|d(\bldr)-d(M(\bldr))|
\leq c q^{-C(i,\bldr,M)}.
\ee
%%%%%%%%%%%%%%%%%%%%%%%%%%%%%%%%%%%%%%%%%%%%%%%%%%%%%

Next, let us look at the growth restrictions coming from the
boundedness of commutators.
In this case, one has the boundedness of only the operators
$[D,\pi(u_{ij})]$. Which means, in effect, one will now have
the condition~(\ref{eqbdd4}) only for $i=1$ and $\bldr\in\Gamma_0$:
\be\label{eqbdd_sph1}
|d(\bldr)-d(M(\bldr))|\leq c q^{-C(1,\bldr,M)}.
\ee
Observe that only allowed moves here are the moves
$M=M_{1,1}\equiv(1)$ and $M=M_{\ell+1,1}\equiv(\ell+1)$.
Looking at the corresponding quantity
$C(1,\bldr,M)$, we find that there are two conditions:
\bea
|d(\bldr^{nk})-d(\bldr^{n,k-1})| &\leq & c,\label{eqbdd_sph2}\\
|d(\bldr^{nk})-d(\bldr^{n+1,k})| &\leq &
       cq^{-\sum_{j=1}^{\ell}H_{1j}(\bldr^{nk})}
      =cq^{-k}.\label{eqbdd_sph3}
\eea
We can now
form a graph by taking $\Gamma_0$ to be the set of vertices,
and by joining two vertices $\bldr$ and $\blds$ by an edge if
$|d(\bldr)-d(\blds)|\leq c$.

\blmma
Let $\scrf_n=\{\bldr^{n,k}:k\in\bbn\}$, $n\in\bbn$.
Then any two points in $\scrf_n$ are connected
by a path lying entirely in $\scrf_n$.

If $n<n'$, then any point in $\scrf_n$ is connected to
any point in $\scrf_{n'}$ by a path such that
$n\leq V_{1,1}(\bldr) \leq n'$
for every vertex $\bldr$ lying on that path.
\elmma
\prf
Take two points $\bldr^{n,j}$ and $\bldr^{n,k}$
in $\scrf_n$. Assume $j<k$.
From the condition (\ref{eqbdd_sph2}), it follows
that any point $\bldr$ is connected to $M_{\ell+1,1}(\bldr)$
by an edge. Therefore the first conclusion follows
from the observation that if we start at $\bldr^{n,k}$
and apply the move $M_{\ell+1,1}$ successively $k-j$ number of times,
we reach the point $\bldr^{n,j}$, and the vertices on this path
are the points $\bldr^{n,i}$ for $i=j, j+1,\ldots,k$.
Observe also that throughout this path, $V_{1,1}(\bldr)$
remains $n$.

For the second part, take a point $\bldr^{n,k}$ in $\scrf_n$
and a point $\bldr^{n',j}$ in $\scrf_{n'}$.
From what we have done above, there is a path
from $\bldr^{n,k}$ to $\bldr^{n,0}$ throughout which
$V_{1,1}(\bldr)=n$.
Similarly there is a path
from $\bldr^{n',j}$ to $\bldr^{n',0}$ throughout which
$V_{1,1}(\bldr)=n'$.
Next, note from (\ref{eqbdd_sph3}) that  for $p\in\bbn$,
the points
$\bldr^{p,0}$ and $\bldr^{p+1,0}$ are connected by an edge
and
$V_{1,1}(\bldr^{p,0})=p$, $V_{1,1}(\bldr^{p+1,0})=p+1$.
So start at $\bldr^{n,0}$ and reach successively the
points
$\bldr^{n+1,0}$, $\bldr^{n+2,0}$ and so on to
eventually reach the point $\bldr^{n',0}$;
also the coordinate $V_{1,1}(\cdot)$ remains between $n$ and $n'$
on this path.\qed

\bthm\label{eqsign_sphere}
Let $D$ be an equivariant Dirac operator on $L_2(G\verb1\1H)$.
Then
\begin{enumerate}
\item
$D$ must be of the form
\[
e_{\bldr,\blds}\mapsto d(\bldr)e_{\bldr,\blds},\quad \bldr\in\Gamma,
\]
where the singular values obey $|d(\bldr)|=O(r_{11})$, and
\item
$\sgn D$ must be of the form $2P-I$
or $I-2P$ where $P$ is, up to a compact perturbation, the projection
onto the closed span of
$\{e_{\bldr^{nk},\blds}: n\in F, k\in\bbn, \blds\in \Gamma_0^{nk}\}$,
for some finite subset $F$ of $\bbn$.
\end{enumerate}
\ethm
\prf
Start with an equivariant self-adjoint operator
$D$ with compact resolvent, so that it is indeed of the form
$e_{\bldr,\blds}\mapsto d(\bldr)e_{\bldr,\blds}$.
By applying a compact perturbation if necessary,
make sure that $d(\bldr)\neq 0$ for all $\bldr\in\Gamma_0$.
We have seen during the proof of the previous lemma that
for any $n$ and $k$ in $\bbn$, the vertices
$\bldr^{nk}$ and $\bldr^{n,k+1}$ are connected by an edge,
and for any $n\in\bbn$, the vertices
$\bldr^{n,0}$ and $\bldr^{n+1,0}$ is connected by an edge.
Thus any vertex $\bldr^{nk}$ can be reached from the vertex
$\bldr^{00}$ by a path of length $n+k$. Therefore one gets the
first assertion.

Next, define
\bean
\Gamma_0^+ &=& \{\bldr\in\Gamma_0: d(\bldr)>0\},\\
\Gamma_0^- &=& \{\bldr\in\Gamma_0: d(\bldr)<0\},\\
\scrf_n^+ &=& \scrf_n\cap \Gamma_0^+,\\
\scrf_n^- &=&  \scrf_n\cap \Gamma_0^-.
\eean
Observe that for the path produced in the proof
of the forgoing lemma to connect two
points $\bldr^{n,k}$ and $\bldr^{n,j}$ in $\scrf_n$,
the coordinate $H_{1,\ell}(\cdot)$ remains between $j$ and $k$.
Now suppose for some  $n$,
both $\scrf_n^+$ and $\scrf_n^-$ are infinite.
Then there are points
\[
0\leq k_1 < k_2 < \ldots
\]
such that $\bldr^{nk}$ is in $\scrf_n^+$ for $k=k_{2j}$
and  $\bldr^{nk}$ is in $\scrf_n^-$ for $k=k_{2j+1}$.
Using the above observation, we can then produce
an infinite ladder by joining
each $\bldr^{n,k_{2j-1}}$ to $\bldr^{n,k_{2j}}$.
Thus for each $n\in\bbn$, exactly one of the sets
$\scrf_n^+$ and $\scrf_n^-$ is finite.
Also, note that by the first part of the previous lemma,
the set of all $n\in\bbn$ for which
both $\scrf_n^+$ and $\scrf_n^-$ are nonempty is finite.
Therefore by applying a compact perturbation, we can ensure that
for every $n$, either $\scrf_n^+=\scrf_n$ or
$\scrf_n^-=\scrf_n$.

Finally, if there are infinitely many $n$'s for which
$\scrf_n^+=\scrf_n$
and infinitely many $n$'s for which $\scrf_n^-=\scrf_n$,
then one can choose a sequence of integers
\[
0\leq n_1 < n_2 <\ldots
\]
such that
$\scrf_n^+=\scrf_n$ for $n=n_{2j}$
and
$\scrf_n^-=\scrf_n$ for $n=n_{2j+1}$.
Now use the second part of the previous lemma
to join each $\bldr^{n_{2j-1},0}$ to $\bldr^{n_{2j},0}$
to produce an infinite ladder.

Thus there is a finite subset $F$ of $\bbn$ such that
exactly one of the following is true:
\[
\scrf_n=\begin{cases}\scrf_n^+ & \mbox{if $n\in F$},\cr
               \scrf_n^- & \mbox{if $n\not\in F$},\end{cases}
\qquad
\mbox{or }
\qquad
\scrf_n=\begin{cases}\scrf_n^- & \mbox{if $n\in F$},\cr
               \scrf_n^+ & \mbox{if $n\not\in F$}.\end{cases}
\]
This is precisely what the second part
of the theorem says.\qed

Next, take the operator $D:e_{\bldr,\blds}\mapsto d(\bldr)e_{\bldr,\blds}$
 on $L_2(G\verb1\1H)$ where the $d(\bldr)$'s are given by:
\be\label{eq_sphere1}
d(\bldr^{nk})=\begin{cases}-k & \mbox{if $n=0$},\cr
                     n+k & \mbox{if $n>0$}.\end{cases}
\ee
\bthm\label{generic_d_sph}
The operator $D$ is an equivariant $(2\ell+1)$-summable
Dirac operator acting on $L_2(G\verb1\1H)$, that gives a
nondegenerate pairing with $K_1(C(G\verb1\1H))$.

The operator $D$ is optimal, i.\ e.\
if $D_0$ is any equivariant Dirac operator on $L_2(G\verb1\1H)$,
then there are positive reals $a$ and $b$ such that
\[
|D_0|\leq a+b|D|.
\]
\ethm
\prf
Recall from equation~(\ref{left_mult}) that the elements
$u_{1,j}$ act on the basis
elements $e_{\bldr^{n,k},\blds}$ as follows:
\bea\label{l2_repn_sph}
u_{1,j}e_{\bldr^{n,k},\blds} &=&
  \sum_{M, M'}
  C_q(1,\bldr^{n,k},M(\bldr^{n,k}))C_q(j,\blds,M'(\blds))
   \kappa(\bldr^{n,k},\blds)
   e_{M(\bldr^{n,k}),M'(\blds)}\cr
 &=& \sum_{M'}
  C_q(1,\bldr^{n,k},M_{11}(\bldr^{n,k}))C_q(j,\blds,M'(\blds))
   \kappa(\bldr^{n,k},\blds)
   e_{M_{11}(\bldr^{n,k}),M'(\blds)} \cr
&&   +
   \sum_{M''}
  C_q(1,\bldr^{n,k},M_{\ell+1,1}(\bldr^{n,k}))C_q(j,\blds,M''(\blds))
   \kappa(\bldr^{n,k},\blds)
   e_{M_{\ell+1,1}(\bldr^{n,k}),M''(\blds)} \cr
&=&\sum_{M'}
  C_q(1,\bldr^{n,k},\bldr^{n+1,k})C_q(j,\blds,M'(\blds))
   \kappa(\bldr^{n,k},\blds)
   e_{\bldr^{n+1,k},M'(\blds)}\cr
&&   +
   \sum_{M''}
  C_q(1,\bldr^{n,k},\bldr^{n,k-1}))C_q(j,\blds,M''(\blds))
   \kappa(\bldr^{n,k},\blds)
   e_{\bldr^{n,k-1},M''(\blds)},
\eea
where the first sum is over all moves $M'\in\bbn^{j}$
whose first coordinate is 1 and the second
sum is over all moves $M''\in\bbn^{j}$
whose first coordinate is $\ell+1$.
If we now plug in the values of the Clebsch-Gordon coefficients
from equations~(\ref{cgc4}) and~(\ref{cgc5}), we get
\bea
u_{1,j}e_{\bldr^{n,k},\blds} &=&
 \sum_{M'}
  P'_1 P'_2 q^{k+C(j,\blds,M')}
   \kappa(\bldr^{n,k},\blds)
   e_{\bldr^{n+1,k},M'(\blds)}\cr
&&     +
   \sum_{M''}
  P''_1 P''_2 q^{C(j,\blds,M'')}
   \kappa(\bldr^{n,k},\blds)
   e_{\bldr^{n,k-1},M''(\blds)},
\eea
where $P'_i$, $P''_j$ and $k(\bldr^{n,k},\blds)$
all lie between two fixed positive numbers.
Boundedness of the commutators $[D,u_{1,j}]$
now follow directly.

For summability, notice that the eigenspace
of $|D|$ corresponding to the eigenvalue $n\in\bbn$
is the span of
\[
\{e_{\bldr^{k,n-k},\blds}: 0\leq k\leq n, \blds\in\Gamma_0^{k,n-k}\}.
\]
Now just count the number of elements in the above set
to get summability.

Next, we will compute the pairing of the $K$-homology class
of this $D$ with a generator of the $K_1$ group.
Write $\omega_q:=q^{-\ell}u_{1,\ell+1}$.
From the commutation relations, it follows that
this element has spectrum
\[
\{z\in\bbc: |z|=0 \mbox{ or }q^n\mbox{ for some }n\in\bbn\}.
\]
Then the element $\gamma_q:=\chi_{\{1\}}(\omega_q^*\omega_q)(\omega_q-I)+I$
is unitary.
We will show that the index of the operator
$Q\gamma_q Q$ (viewed as an operator on $QL_2(G\verb1\1H)$)
is $1$, where $Q=\frac{I-\sgn D}{2}$, i.\ e.\ it
is the projection onto the closed
linear span of
$\{e_{\bldr^{0,k},\blds}:k\in\bbn, \blds\in\Gamma_0^{0,k}\}$.
What we will actually do is
compute the index of the
operator $Q\gamma_0 Q$ and appeal to continuity of the index.
From equation~(\ref{l2_repn_sph}), we get
\bea\label{for_q=0_1}
\lefteqn{u_{1,\ell+1}e_{\bldr^{0,k},\blds}}\cr
 &=&
  C_q(1,\bldr^{0,k},M_{11}(\bldr^{0,k}))C_q(\ell+1,\blds,N_{1,0}(\blds))
    \kappa(\bldr^{0,k},M_{11}(\bldr^{0,k}))e_{\bldr^{1,k},N_{1,0}(\blds)}\cr
&& +     C_q(1,\bldr^{0,k},M_{\ell+1,1}(\bldr^{0,k}))
                C_q(\ell+1,\blds,M_{\ell+1,\ell+1}(\blds))
    \kappa(\bldr^{0,k},M_{\ell+1,1}(\bldr^{0,k}))
        e_{\bldr^{0,k-1},M_{\ell+1,\ell+1}(\blds)}.\cr
&&
\eea
Use the formula~(\ref{cgc1}) for Clebsch-Gordon coefficients
to get
\bea
C_q(1,\bldr^{0,k},M_{11}(\bldr^{0,k}))
 &=& q^{k}(1+ o(q)),\\
C_q(1,\bldr^{0,k},M_{\ell+1,1}(\bldr^{0,k})
 &=&  1+ o(q),\\
C_q(\ell+1,\blds,N_{1,0}(\blds))
 &=& 1+ o(q),\\
C_q(\ell+1,\blds,M_{\ell+1,\ell+1}(\blds))
 &=& q^{s_{\ell+1,1}+\ell}(1+ o(q)),
\eea
where $o(q)$ signifies a function of $q$ that is
continuous at $q=0$ and
$o(0)=0$.
We also have
\bea
\kappa(\bldr^{0,k},M_{11}(\bldr^{0,k}))
  &=& q^\ell (1+o(q)),\\
\kappa(\bldr^{0,k},M_{\ell+1,1}(\bldr^{0,k}))
  &=& 1+o(q),
\eea
where $o(q)$ is as earlier.
Plugging these values in~(\ref{for_q=0_1}) 
we get
\be
\omega_q e_{\bldr^{0,k},\blds}
  =  q^{k}(1+o(q))e_{\bldr^{1,k},N_{1,0}(\blds)}
   + q^{s_{\ell+1,1}}(1+o(q))e_{\bldr^{0,k-1},M_{\ell+1,\ell+1}(\blds)}
\ee
Putting $q=0$,
we get
\be
\omega_0 e_{\bldr^{0,k},\blds}
  = \begin{cases}
   e_{\bldr^{0,k-1},M_{\ell+1,\ell+1}(\blds)}
             & \mbox{if $k>0$ and $s_{\ell+1,1}=0$},\cr
    e_{\bldr^{1,0},N_{1,0}(\blds)} & \mbox{if $k=0$},\cr
    0 & \mbox{otherwise.}
  \end{cases}
\ee
Thus $\omega_0^*\omega_0$ is the projection onto the span
of
$\{e_{\bldr^{0,k},\blds^k}: k\in\bbn\}$
where
$\blds^k$ is the GT tableaux given by
\[
s^k_{ij}=\begin{cases}
             0 & \mbox{if $i=\ell+2-j$},\cr
                  k & \mbox{otherwise},\end{cases}
\]
which is uniquely determined by the conditions
$s_{\ell+1,1}=0$ and that $\blds\in\Gamma_0^{0,k}$.
Therefore the operator $\gamma_0$  is given by
\[
\gamma_0 e_{\bldr^{0,k},\blds} =
 e_{\bldr^{0,k},\blds} - \chi_{\{\blds=\blds^k\}}e_{\bldr^{0,k},\blds}
   + \chi_{\{\blds=\blds^k\}}
                   e_{\bldr^{0,k-1},\blds^{k-1}}.
\]
It now follows that the index of $Q\gamma_0 Q$ is $1$.

Note that if $D_0$ is an equivariant Dirac operator
with eigenvalues $d_0(\bldr)$, then by theorem~\ref{eqsign_sphere} there is a
$b>0$ such that
\[
 |d_0(\bldr)|<b r_{11}=b|d(\bldr)|,\quad \bldr\neq \mathbf{0}.
\]
Write $a=|d_0(\mathbf{0})|$. Then we have the required inequality.
\qed

\vspace{2ex}

\noindent
{\footnotesize \textbf{Acknowledgement}: The second author would like to thank the
Isaac Newton Institute and in particular the organisers of the NCG semester for their hospitality.}

%%%%%%%%%%%%%%%%%%%%%%%%%%
%%%  BIBLIOGRAPHY
%%%%%%%%%%%%%%%%%%%%%%%%%%

%%%%%%%%%%%%%%%%%%%%%%%%%%%%%%%
\noindent{\sc Partha Sarathi Chakraborty}
(\texttt{parthac@imsc.res.in})\\
         {\footnotesize  Institute of Mathematical Sciences, 
CIT Campus, Chennai--600\,113, INDIA}\\[1ex]
{\sc Arupkumar Pal} (\texttt{arup@isid.ac.in})\\
         {\footnotesize Indian Statistical
Institute, 7, SJSS Marg, New Delhi--110\,016, INDIA}

%%%%%%%%%%%%%%%%%%%%%%%%%%%%%%%%%%%%%%%%%%%%%%%%%%%%%%%%%%%
\end{document}